\documentclass[10pt]{article}
\usepackage{amsfonts}
\def\frak{\fontencoding{U}\fontfamily{euf}\selectfont}
\begin {document}
\topmargin= -.2in \baselineskip=16pt

\title {Trivial Factors For $L$-functions of Symmetric Products
of Kloosterman Sheaves}

\author {Lei Fu\\
{\small Chern Institute of Mathematics and LPMC, Nankai University,
Tianjin, P.
R. China}\\
{\small leifu@nankai.edu.cn}\\{}\\
Daqing Wan\\
{\small Department of Mathematics, University of California,
Irvine, CA
92697}\\
{\small dwan@math.uci.edu}}
\date{}
\maketitle

\bigskip
\centerline {\bf 0. Introduction}

\bigskip
\bigskip
In this paper, we determine the trivial factors of $L$-functions of
both finite and infinite symmetric products of Kloosterman sheaves.

Let ${\bf F}_q$ be a finite field of characteristic $p$ with $q$
elements, let $l$ be a prime number distinct from $p$, and let
$\psi: {\bf F}_q\to \overline {\bf Q}_l^\ast$ be a nontrivial
additive character. Fix an algebraic closure ${\bf F}$ of ${\bf
F}_q$. For any integer $k$, let ${\bf F}_{q^k}$ be the extension
of ${\bf F}_q$ in ${\bf F}$ with degree $k$. Let $n\geq 2$ be a
positive integer. If $\lambda$ lies in ${\bf F}_{q^k}$, we define
the $(n-1)$-variable Kloosterman sum by
$${\rm Kl}_{n}({\bf F}_{q^k}, \lambda)=\sum\limits_{x_1\cdots x_{n}=\lambda
,\;x_i\in  {\bf F}_{q^k}} \psi({\rm Tr}_{{\bf F}_{q^k}/{\bf
F}_q}(x_1+\cdots + x_{n} )).$$ Such character sums can be studied
via either $p$-adic methods or $l$-adic methods. In [D1]
Th\'eor\`eme 7.8, Deligne constructs a lisse $\overline{\bf
Q}_l$-sheaf of rank $n$ on ${\bf A}^1_{{\bf F}_q}-\{0\}$ pure of
weight $n-1$, which we denote by ${\rm Kl}_{n}$ and call the
Kloosterman sheaf, with the property that for any $x\in ({\bf
A}^1_{{\bf F}_q}-\{0\})({\bf F}_{q^k})={\bf F}_{q^k}^\ast$, we have
$${\rm Tr}(F_x, {\rm Kl}_{n,\bar x})=(-1)^{n-1} {\rm Kl}_{n}({\bf F}_{q^k},
x),$$ where $F_x$ is the geometric Frobenius element at the point
$x$. Let $\eta$ be the generic point of ${\bf A}^1_{{\bf F}_q}$.
The Kloosterman sheaf gives rise to a Galois representation
$${\rm Kl}_{n}:{\rm Gal}(\overline {{\bf F}_q(T)}/{\bf F}_q(T))\to
{\rm GL}(({\rm Kl}_n)_{\bar\eta})$$ unramified outside $0$ and
$\infty$. From the $p$-adic point of view, the Kloosterman sheaf is
given by an ordinary overconvergent F-crystal of rank $n$ over ${\bf
A}^1_{{\bf F}_q}-\{0\}$. See Sperber [S].

For each positive integer $k$, denote the $L$-function of the $k$-th
symmetric product of the Kloosterman sheaf by $L(k,n,T)$:
$$L(k,n,T): =L({\bf A}_{{\bf F}_q}^1-\{0\},
{\rm Sym}^k{\rm Kl}_n,T) \in 1+T{\bf Z}[[T]],$$ and we call it
simply the $k$-th symmetric product $L$-function. This is a rational
function whose reciprocal zeros and poles are Weil $q$-numbers by
theorems of Grothendieck and Deligne. Our aim of this paper is to
understand the trivial factors of $L(k,n,T)$ and their variation as
the integer $k$ varies $p$-adically. In the case $n=2$, the trivial
factor problem for $L(k, 2, T)$ was first studied by Robba ([R]) via
Dwork's $p$-adic cohomology. Robba determined the trivial factors
for $L(k, 2, T)$ assuming $p>k/2$. Using $l$-adic results of Deligne
and Katz, we ([FW]) determined the trivial factor at $\infty$ for
general $L(k,n,T)$. The trivial factor at $0$ is easy to determine
for $n=2$. But in [FW] we were unable to determine the trivial
factor at $0$ of $L(k,n,T)$ for $n>2$. We solve this problem in the
present paper.

In a different but related direction, the $p$-adic limit of
$L(k,n,T)$ when $k$ goes to infinity in a fixed $p$-adic direction
was shown to be a $p$-adic meromorphic function in [W1]. This idea
was the key in proving Dwork's unit root conjecture for the
Kloosterman family. See [W1], [W2] and [W3]. To be precise, for a
$p$-adic integer $s$, we choose a sequence of positive integers
$k_i$ which approaches $s$ as $p$-adic integers but goes to infinity
as complex numbers. Then we define the $p$-adic $s$-th symmetric
product $L$-function to be
$$L_p(s,n,T) =\lim_{i\rightarrow \infty} L(k_i,n,T)\in 1+T{\bf Z}_p[[T]].$$
This limit exists as a formal $p$-adic power series and is
independent of the choice of the sequence $k_i$. It is a sort of two
variable $p$-adic $L$-function.  Note that even when $s$ is a
positive integer, $L_p(s,n,T)$ is very different from $L(s,n,T)$. It
was shown in [W1] that $L_p(s, n, T)$ is a $p$-adic meromorphic
function by a uniform limiting argument. Alternatively, it was shown
in [W2] that
$$L_p(s,n,T) = L(M_s(\infty), T),$$ where $M_s(\infty)$ is an infinite rank
nuclear overconvergent $\sigma$-module on ${\bf A}_{{\bf
F}_q}^1-\{0\}$. This gives another proof that $L_p(s,n,T)$ is
$p$-adic meromorphic.

In this paper, combining $p$-adic methods and $l$-adic methods, we
prove the following more precise result.

\bigskip
\noindent {\bf Theorem 0.1.} For each $p$-adic integer $s$, the
$p$-adic $s$-th symmetric product $L$-function $L_p(s,n,T)$ is a
$p$-adically entire function (i.e., no poles). Furthermore, the
entire function $L_p(s,n,T)$ is divisible by the $p$-adic entire
function $\prod_{i=0}^{\infty} (1-q^iT)^{d_i}$, where $d_j$ is the
coefficient of $x^j$ in the power series expansion of
$${1 \over (1-x^2)(1-x^3)\cdots (1-x^{n-1})},$$
that is, for each $s\in {\bf Z}_p$, the entire function $L_p(s,n,T)$
has a zero at $T= q^{-j}$ with multiplicity at least $d_j$ for each
non-negative integer $j$.

\bigskip
\noindent {\bf Remark}. Grosse-Kl\"onne [GK] showed the $p$-adic
meromorphic continuation of $L_p(s,n, T)$ to some $s\in {\bf Q}_p$
with $|s|_p < 1+\epsilon$ for some small $\epsilon >0$. We do not
know if Theorem 0.1 can be extended to such non-integral $p$-adic
$s$.
\bigskip

In order to prove the above theorem for infinite symmetric product
$L$-function $L_p(s,n,T)$, we need to have a good understanding of
the finite symmetric product $L$-function $L(k,n,T)$ for every
positive integer $k$. Let $j:{\bf A}_{{\bf F}_q}^1-\{0\}\to {\bf
P}^1_{{\bf F}_q}$ be the canonical open immersion. By definition, we
have the following relation between the $L$-functions $L(k,n,T)$ and
$L({\bf P}_{{\bf F}_q}^1, j_\ast({\rm Sym}^k({\rm Kl}_n)),T)$:
\begin{eqnarray*}
&&L(k,n,T)\\
&=&L({\bf P}_{{\bf F}_q}^1, j_\ast({\rm Sym}^k({\rm Kl}_n)),T){\rm
det}(1-F_0T,({\rm Sym}^k({\rm Kl}_n)_{\bar \eta})^{I_0}){\rm
det}(1-F_\infty T,({\rm Sym}^k({\rm Kl}_n)_{\bar \eta})^{I_\infty}),
\end{eqnarray*}
where $I_0$ (resp. $I_\infty$) is the inertia subgroup at $0$ (resp.
$\infty$), and $F_0$ (resp. $F_\infty$) is the geometric Frobenius
element at $0$ (resp. $\infty$). Here we use the fact that
$$({\rm Sym}^k({\rm
Kl}_n)_{\bar \eta})^{I_0}=(j_\ast ({\rm Sym}^k({\rm Kl}_n)))_{\bar
0}, \; ({\rm Sym}^k({\rm Kl}_n)_{\bar \eta})^{I_\infty}=(j_\ast
({\rm Sym}^k({\rm Kl}_n)))_{\overline \infty}.$$ We call ${\rm
det}(1-F_0T,({\rm Sym}^k({\rm Kl}_n)_{\bar \eta})^{I_0})$ (resp.
${\rm det}(1-F_\infty T,({\rm Sym}^k({\rm Kl}_n)_{\bar
\eta})^{I_\infty})$) the local factor at $0$ (resp. $\infty$) of
$L(k,n,T)$. On the other hand, by Grothendieck's formula for
$L$-functions, we have
\begin{eqnarray*}
&&L({\bf P}_{{\bf F}_q}^1, j_\ast({\rm Sym}^k({\rm
Kl}_n)),T)\\
&=&\frac{{\rm det}(1-FT, H^1({\bf P}_{\bf F}^1,j_\ast({\rm
Sym}^k({\rm Kl}_n))))}{{\rm det}(1-FT, H^0({\bf P}_{\bf
F}^1,j_\ast({\rm Sym}^k({\rm Kl}_n)))){\rm det}(1-FT, H^2({\bf
P}_{\bf F}^1,j_\ast({\rm Sym}^k({\rm Kl}_n))))}.
\end{eqnarray*}
So we get the factorization
\begin{eqnarray*}
&&L(k,n, T)\\&=& \frac{{\rm det}(1-FT, H^1({\bf P}_{\bf
F}^1,j_\ast({\rm Sym}^k({\rm Kl}_n)))){\rm det}(1-F_0T,(({\rm
Sym}^k{\rm Kl}_n)_{\bar \eta})^{I_0}){\rm det}(1-F_\infty T,(({\rm
Sym}^k{\rm Kl}_n)_{\bar \eta})^{I_\infty})}{{\rm det}(1-FT, H^0({\bf
P}_{\bf F}^1,j_\ast({\rm Sym}^k({\rm Kl}_n)))){\rm det}(1-FT,
H^2({\bf P}_{\bf F}^1,j_\ast({\rm Sym}^k({\rm Kl}_n))))}.
\end{eqnarray*}
The first factor ${\rm det}(1-FT, H^1({\bf P}_{\bf F}^1,j_\ast({\rm
Sym}^k({\rm Kl}_n))))$ is called the non-trivial factor. It is pure
of weight $k(n-1)+1$ by [D2] 3.3.1. All other factors on the right
side of the above expression are called trivial factors. The zeros
of these trivial factors give rise to the trivial zeros or poles of
$L(k,n,T)$. We will determine all the trivial factors and their
variation with $k$ as $k$ varies. As a consequence, it gives some
partial information on the non-trivial factor and its variation with
$k$ as well.

We now describe the results for the trivial factors. In [FW], we
studied in detail the behavior of the Kloosterman representation
at $\infty$, and we used our results to calculate the local factor
at $\infty$ of $L(k,n,T)$. To recall our result, let $\zeta$ be a
primitive $n$-th root of unity in ${\bf F}$. For each positive
integer $k$, let $S_k(n,p)$ be the set of $n$-tuples $(j_0,
\cdots, j_{n-1})$ of non-negative integers satisfying
$j_0+j_1+\cdots +j_{n-1}=k$ and $j_0+j_1\zeta +\cdots
+j_{n-1}\zeta^{n-1}=0$ in ${\bf F}$. Let $\sigma$ be the cyclic
shifting operator
$$\sigma(j_0, \cdots, j_{n-1}) = (j_{n-1}, j_0, \cdots, j_{n-2}).$$
It is clear that the set $S_k(n,p)$ is $\sigma$-stable. Let $V$ be a
$\overline{\bf Q}_l$-vector space of dimension $n$ with basis
$\{e_0, \cdots, e_{n-1}\}$. For an $n$-tuple $j=(j_0,\cdots,
j_{n-1})$ of non-negative integers such that $j_0+\cdots
+j_{n-1}=k$, write
$$e^j =e_0^{j_0}e_1^{j_1}\cdots e_{n}^{j_{n-1}}$$
as an element of ${\rm Sym}^kV$. For such an $n$-tuple $j$, we
define
$$v_j = \sum_{i=0}^{n-1} (-1)^{j_{n-1}+\cdots +j_{n-i}} e^{\sigma^i(j)}.$$
If $k=j_0+j_1+\cdots+ j_{n-1}$ is even, then we have
$v_{\sigma(j)}=(-1)^{j_{n-1}}v_j$. Let $a_k(n,p)$ be the number of
$\sigma$-orbits in $S_k(n,p)$. When $k$ is even, let $b_k(n,p)$ be
the number of those $\sigma$-orbits $j$ in $S_k(n,p)$ such that
$v_j\not=0$, and let $c_k(n,p)$ be the number of $\sigma$-orbits $j$
in $S_k(n,p)$ such that $v_j\not=0$ and that $j_1+2j_2+\cdots
+(n-1)j_{n-1}$ is odd. Our result on the local factor at $\infty$ of
$L(k,n,T)$ is the following.

\bigskip
\noindent{\bf Theorem 0.2.} (Theorem 2.5 in [FW]) Suppose
$n|(q-1)$.

(1) If $n$ is odd, then for all $k$, we have
\[{\rm det}(1-F_\infty T, ({\rm Sym}^k({\rm Kl}_{n})_{\bar\eta})^{I_\infty})=
(1-q^{\frac{k(n-1)}{2}}T)^{a_k(n,p)}.\]

(2) If $n$ is even and $k$ is odd, then we have
\[{\rm det}(1-F_\infty T, ({\rm Sym}^k({\rm Kl}_{n})_{\bar\eta})^{I_\infty})= 1.\]

(3) Suppose $n$ and $k$ are both even. We have
\begin{eqnarray*}
&&{\rm det}(1-F_\infty T, ({\rm Sym}^k({\rm Kl}_{n})_{\bar\eta})^{I_\infty})\\
&=& \left\{
\begin {array} {ll}
(1-q^{\frac{k(n-1)}{2}}T)^{b_k(n,p)}& {\rm if} \; 2n|(q-1), \\
(1+q^{\frac{k(n-1)}{2}}T)^{c_k(n,p)}(1-q^{\frac{k(n-1)}{2}}T)^
{b_k(n,p)-c_k(n,p)}& {\rm if} \; 2n\not|(q-1), \hbox { either }
4|n\;
{\rm  or}\; 4|k,\\
(1-q^{\frac{k(n-1)}{2}}T)^{c_k(n,p)}(1+q^{\frac{k(n-1)}{2}}T)^
{b_k(n,p)-c_k(n,p)}& {\rm if} \; 2n\not|(q-1), \;4\not|n \;{\rm
and}\; 4\not|k.
\end{array}\right.
\end{eqnarray*}

\bigskip
In this paper, we get the following formula for the local factor at
$0$ of $L(k,n,T)$.

\bigskip
\noindent {\bf Theorem 0.3.} We have
$${\rm det}(I-F_0T, ({\rm Sym}^k({\rm Kl}_{n})_{\bar\eta})^{I_0})=\prod_{u=0}^
{[\frac{k(n-1)}{2}]} (1-q^uT)^{m_k(u)},$$ where $m_k(u)$ is
determined by
$$\frac{(1-x^{n})\cdots (1-x^{n+k-2})(1-x^{n+k-1})}{(1-x^2)\cdots
(1-x^{k-1})(1-x^k)}=\sum_{u=0}^\infty m_k(u)x^u.$$ We have
$$m_k(u)=c_k(u)-c_k(u-1),$$ where $c_k(u)$ is the number of elements of the set
$$\{(i_0,\ldots, i_{n-1})| i_j\geq 0, \;
i_0+i_1+\cdots +i_{n-1}=k,\;0\cdot i_0+1\cdot i_1+\cdots +(n-1)\cdot
i_{n-1}=u\}.$$

\bigskip
The trivial poles of $L(k,n,T)$ can be derived from Katz's global
monodromy theorem and Grothendieck's formula for $L$-functions. For
completeness, we include this deduction in detail by working out the
relevant representation theory which should be well known to
experts.

Denote by $G$ the Zariski closure of the image of ${\rm
Gal}(\overline {{\bf F}(T)}/{\bf F}(T))$ under the representation
$${\rm Kl}_{n}:{\rm Gal}(\overline {{\bf F}_q(T)}/{\bf F}_q(T))\to
{\rm GL}(({\rm Kl}_n)_{\bar\eta}).$$ By [K] 11.1, we have
$$
G=\left\{\begin{array}{cl} {\rm Sp}(n)&\hbox { if } n \hbox { is even},\\
{\rm SL}(n)&\hbox { if } n \hbox { is odd},\hbox { and } p\not=2.
\\{\rm SO}(n)&\hbox { if } n \hbox { is odd, } n\not= 7 \hbox { and } p=2,
\\{\rm G}_2&\hbox { if } n=7 \hbox { and } p=2.
\end{array}\right.$$
If $pn$ is even, we have $(-1)^n=1$ in ${\bf F}_q$. By [K] 4.2.1, we
then have a perfect paring
$${\rm Kl}_n\otimes {\rm Kl}_n\to \overline {\bf Q}_l(1-n).$$
When $n$ is even, we have $G={\rm Sp}(n)$, the paring is
alternating, and ${\rm Kl}_n$ is isomorphic to the standard
representation of ${\rm Sp}(n)$. When $n$ is odd and $p=2$, we have
$G={\rm SO}(n)$ or $G={\rm G}_2$, the paring is symmetric, and ${\rm
Kl}_n$ is isomorphic to the standard representation of the ${\rm
SO}(n)$ or ${\rm G}_2$. (The standard representation of ${\rm G}_2$
is defined to be the unique irreducible representation of dimension
7.) When $pn$ is odd, ${\rm Kl}_n$ is isomorphic to the standard
representation of ${\rm SL}(n)$. In the Appendix of this paper, we
will prove the following result:

\bigskip
\noindent {\bf Lemma 0.4.} Let $\hbox {\frak g}$ be one of the
following Lie algebras
$$\hbox{\frak sl}(n), \hbox {\frak sp}(n),
\hbox{\frak so}(n), \hbox{\frak g}_2$$ and let $V$ be the standard
representation of $\hbox{\frak g}$. In the case where $\hbox{\frak
g}=\hbox {\frak sl}(n)$ or $\hbox {\frak sp}(n)$, the representation
${\rm Sym}^kV$ is irreducible, and in the case where $\hbox{\frak
g}=\hbox{\frak so}(n)$ or $\hbox{\frak g}_2$, the representation
${\rm Sym}^kV$ contains exactly one copy of the trivial
representation if $k$ is even, and contains no trivial
representation if $k$ is odd.

\bigskip
By [D2] 1.4.1, we have
\begin{eqnarray*}
&& H^0({\bf P}_{\bf F}^1,j_\ast({\rm Sym}^k({\rm Kl}_n)))=(({\rm
Kl}_n)_{\bar \eta})^{{\rm Gal}(\overline{{\bf F}_q(t)}/{\bf
F}_q(t))}=(({\rm Kl}_n)_{\bar \eta})^G,\\
&& H^2({\bf P}_{\bf F}^1,j_\ast({\rm Sym}^k({\rm Kl}_n)))=(({\rm
Kl}_n)_{\bar \eta})_{{\rm Gal}(\overline{{\bf F}_q(t)}/{\bf
F}_q(t))}(-1)=(({\rm Kl}_n)_{\bar \eta})_G(-1).
\end{eqnarray*}
Combined with [K] 11.1 and Lemma 0.4, we get the following.

\bigskip
\noindent {\bf Corollary 0.5.} We have
\begin{eqnarray*}
{\rm det}(1-FT, H^0({\bf P}_{{\bf F}}^1,j_\ast ({\rm Sym}^k({\rm
Kl}_n))))&=&\left\{\begin{array}{cl} 1 & \hbox { if } n \hbox { is
even, or
} k \hbox { is odd, or } pn \hbox { is odd},\\
1-q^{\frac{k(n-1)}{2}}T&\hbox { if } p=2, \hbox { $k$ is even and
} n \hbox { is odd},
\end{array}\right. \\
{\rm det}(1-FT, H^2({\bf P}_{{\bf F}}^1,j_\ast( {\rm Sym}^k({\rm
Kl}_n))))&=&\left\{\begin{array}{cl} 1 & \hbox { if } n \hbox { is
even, or
} k \hbox { is odd, or } pn \hbox { is odd},\\
1-q^{\frac{k(n-1)+2}{2}}T&\hbox { if } p=2, \hbox { $k$ is even
and } n \hbox { is odd}.
\end{array}\right.
\end{eqnarray*}

\bigskip
Combining these results with the $p$-adic limiting argument in [W1],
we shall obtain infinitely many trivial zeros (if $n>2$) for the
infinite $p$-adic symmetric product $L$-function $L_p(s,n,T)$ as
stated in Theorem 0.1. This suggests that there should be an
interesting trivial zero theory for the $L$-function of any infinite
$p$-adic symmetric product of a pure $l$-adic sheaf whose $p$-adic
unit root part has rank one. Our result here provides the first
evidence for such a theory.

\bigskip
The paper is organized as follows. In \S 1, we recall the canonical
form of the local monodromy of the Kloosterman sheaf at $0$. In \S
2, we summarize the basic representation theory for $\hbox{\frak sl}
(2)$. In \S 3, we prove Theorem 0.3 using results in the previous
two sections. In \S 4, we apply our results on local factors at $0$
to prove Theorem 0.1. In section 5, we derive some consequences for
the non-trivial factors and its variation with $k$. In the appendix,
we include a proof of Lemma 0.4 which implies Corollary 0.5.

\bigskip\bigskip
\noindent {\bf Acknowledgements.} The research of Lei Fu is
supported by NSFC (10525107). The research of Daqing Wan is
partially supported by NSF.

\bigskip
\bigskip
\centerline {\bf 1. The Canonical Form of the Local
Monodromy}

\bigskip
\bigskip
Let $K$ be a local field with residue field ${\bf F}_q$, and let
$$\rho: {\rm Gal}(\overline K/K)\to {\rm GL}(V)$$ be a
$\overline {\bf Q}_l$-representation. Suppose the inertia subgroup
$I$ of ${\rm Gal}(\overline K/K)$ acts unipotently on $V$. Fix a
uniformizer $\pi$ of $K$, and consider the $l$-adic part of the
cyclotomic character
$$t_l: I\to {\bf Z}_l(1), \; \sigma\mapsto
\left(\frac{\sigma(\sqrt[l^n]{\pi})}{\sqrt[l^n]{\pi}}\right).$$
Note that for $\sigma$ in the inertia group, the $l^n$-th root of
unity $\frac{\sigma(\sqrt[l^n]{\pi})}{\sqrt[l^n]{\pi}}$ does not
depends on the choice of the $l^n$-th root $\sqrt[l^n]{\pi}$ of
$\pi$. Since the restriction to $I$ is unipotent, there exists a
nilpotent homomorphism
$$N:V(1)\to V$$ such that
$$\rho(\sigma)=\exp (t_l(\sigma).N)$$ for any $\sigma\in I$.
Fix a lifting $F\in {\rm Gal}(\overline K/K)$ of the geometric
Frobenius element in ${\rm Gal}({\bf F}/{\bf F}_q)$. We have
$$t_l(F^{-1}\sigma F)=t_l(\sigma)^q.$$ So
\begin{eqnarray*}
\exp(t_l(\sigma).N)\rho(F)&=&
\rho(\sigma)\rho(F)\\
&=&\rho(\sigma F)\\
&=&\rho (F F^{-1}\sigma F)\\
&=&\rho(F)\rho(F^{-1}\sigma F)\\
&=&\rho(F)\exp(t_l(F^{-1}\sigma F).N)\\
&=&\rho(F)\exp(qt_l(\sigma).N).
\end{eqnarray*}
Therefore
$$\rho(F)^{-1}\exp(t_l(\sigma).N)\rho(F)=\exp(qt_l(\sigma).N).$$
Hence
$$\rho(F)^{-1}(t_l(\sigma).N)\rho(F)=qt_l(\sigma).N.$$
Fix a generator $\zeta$ of ${\bf Z}_l(1)$. Choose $\sigma\in I$ so
that $t_l(\sigma)=\zeta$. For convenience, denote $\rho(F)$ by
$F$, and denote the homomorphism $$V\to V, v\mapsto
N(v\otimes\zeta)$$ by $N$. Then the last equation gives
$$F^{-1}NF=qN,$$
that is,
$$NF=qFN.$$

\bigskip
Now let's take $K$ to be the completion of ${\bf F}_q(T)$ at $0$,
and let $\rho:{\rm Gal}(\overline K/K)\to {\rm GL}(V)$ be the
restriction of the representation ${\rm Kl}_{n}:{\rm Gal}(\overline
{{\bf F}_q(T)}/{\bf F}_q(T))\to {\rm GL}(({\rm Kl}_n)_{\bar \eta})$
defined by the Kloosterman sheaf. In [D1] Th\'eor\`eme 7.8, it is
shown that the inertia group $I_0$ at $0$ acts unipotently on $({\rm
Kl}_n)_{\bar\eta}$ with a single Jordan block, and the geometric
Frobenius $F_0$ at $0$ acts trivially on the invariant $(({\rm
Kl}_n)_{\bar\eta})^{I_0}$ of the inertia group. With the above
notations, this means the nilpotent map $N$ has a single Jordan
block, and $F$ acts trivially on ${\rm ker}(N)$. By [D2] 1.6.14.2
and 1.6.14.3, the eigenvalues of $F$ are $1,q,\ldots, q^{n-1}$.
Since the number of distinct eigenvalues is exactly the rank of the
Kloosterman sheaf, $F$ is diagonalizable. Let $v$ be a (nonzero)
eigenvector of $F$ with eigenvalue $q^{n-1}$. Using the equation
$NF=qFN$, we see $N(v)$ is an eigenvector of $F$ with eigenvalue
$q^{n-2}$. Note that if $n\geq 2$, then $N(v)$ can not be $0$ since
otherwise $v$ lies in ${\rm ker}(N)$ and $F$ does not acts trivially
on $v$. This contradicts to the property of the Kloosterman sheaf.
Similarly, if $n\geq 3$, then $N^2(v)$ is a nonzero eigenvector of
$F$ with eigenvalue $q^{n-3}$, $\ldots$, and $N^{n-1}(v)$ is a
nonzero eigenvector of $F$ with eigenvalue $1$, and $N^{n}(v)=0$. As
$v, N(v), \ldots, N^{n-1}(v)$ are nonzero eigenvectors of $F$ with
distinct eigenvalues, they are linearly independent and form a basis
of $V$. With respect to the basis $\{N^{n-1}(v),\ldots, N(v), v\}$,
the matrix of $N$ is
$$\left(\begin{array}{cccc}
0&1&&\\
&0&\ddots&\\
&&\ddots&1\\
&&&0
\end{array}
\right),$$ and the matrix of $F$ is
$$\left(\begin{array}{cccc}
1&&&\\
&q&&\\
&&\ddots&\\
&&&q^{n-1}
\end{array}
\right).$$ We summarize the above results as follows.

\bigskip
\noindent {\bf Proposition 1.1.} Notation as above. For the triple
$(V, F, N)$ defined by the Kloosterman sheaf, there exists a basis
$e_0,\ldots, e_{n-1}$ of $V$ such that
$$F(e_0)=e_0,\; F(e_1)=qe_1,\; \ldots, \;F(e_{n-1})=q^{n-1}e_{n-1}$$ and
$$N(e_0)=0,\; N(e_1)=e_0,\; \ldots,\; N(e_{n-1})=e_{n-2}.$$

\bigskip
\bigskip
\centerline {\bf 2. Representation of $\hbox{\frak sl}(2)$}

\bigskip
In this section, we summarize the representation theory of the Lie
algebra $\hbox{\frak sl}(2)$ of traceless matrices over the field
$\overline {\bf Q}_l$. Consider the following three elements
$$H=\left(\begin{array}{cc}
1&0\\
0&-1
\end{array}
\right),\;X=\left(\begin{array}{cc}
0&1\\
0&0
\end{array}
\right),\;Y=\left(\begin{array}{cc}
0&0\\
1&0
\end{array}
\right).$$ We have
$$[X,Y]=H,\; [H,X]=2X,\; [H,Y]=-2Y.$$
Let $V$ be a finite dimensional irreducible $\overline {\bf
Q}_l$-representation of $\hbox{\frak sl}(2)$, and let $u\in V$ be a
(nonzero) eigenvector of the action $H$ on $V$ with eigenvalue
$\lambda$. Using the above relations, we get
$$HXu=(\lambda+2)Xu,\; HYu=(\lambda-2)Yu.$$
It follows that $u, Xu, X^2u, \ldots$ are eigenvectors of $H$ with
different eigenvalues $\lambda, \lambda+2,\lambda+4,\ldots$ Since
$V$ is finite dimensional, we must have $X^ku=0$ for sufficiently
large integer $k$. Let $m$ be the largest integer such that
$X^{m}u\not=0$. It is an eigenvector of $H$ with eigenvalue
$\lambda+2m$. Set $v=X^mu$ and $\mu=\lambda+2m$. Then $Xv=0$ and
$v,Yv,Y^2v,\ldots$ are eigenvectors of $H$ with different
eigenvalues $\mu, \mu-2,\mu-4,\ldots$ Let $n$ be the largest integer
such that $Y^{n}v\not=0$. Then $\{v,Yv,\ldots, Y^nv\}$ are linearly
independent. The space ${\rm span}\{v,Yv,\ldots, Y^nv\}$ is
invariant under the actions of $H$ and $Y$. Moreover, we have
\begin{eqnarray*}
&&Xv=0, \\
&&X(Yv)=Y(Xv)+[X,Y]v=Hv=\mu v,\\
&&X(Y^2v)=Y(X(Yv))+[X,Y](Yv)=(\mu+(\mu-2))Yv,\\
&&\ldots\\
&&X(Y^iv)=Y(X(Y^{i-1}v))+[X,Y](Y^{i-1}v)\\
&&\qquad\quad\;\,=(\mu+(\mu-2)+\cdots+(\mu-2(i-1)))Y^{i-1}v,\\
&&\qquad\quad\;\,= i(\mu-i+1)Y^{i-1}v,\\
&&\ldots
\end{eqnarray*}
So ${\rm span}\{v,Yv,\ldots, Y^nv\}$ is also invariant under the
action of $X$. Hence ${\rm span}\{v, Yv,\ldots, Y^nv\}$ is
invariant under the action of $\hbox{\frak sl}(2)$. Since $V$ is
an irreducible representation of $\hbox{\frak sl}(2)$, it follows
that $\{v,Yv,\ldots, Y^nv\}$ is a basis of $V$. By the choice of
$n$, we have $Y^nv\not=0$ and $Y^{n+1}v=0$. On the other hand, we
have
$$0=X(Y^{n+1}v)=(n+1)(\mu-n)Y^{n}v.$$ So $(n+1)(\mu-n)=0$, and hence $\mu=n$.
We summarize our results as follows.

\bigskip
\noindent {\bf Proposition 2.1.} Let $V$ be a finite dimensional
irreducible $\overline{\bf Q}_l$-representation of $\hbox{\frak
sl}(2)$. Then there exists a (nonzero) eigenvector $v$ of $H$ such
that $Xv=0$. Such a vector is called a {\it highest weight vector}
for the representation $V$. Let $n={\rm dim}(V)-1$. For any highest
weight vector $v$, we have
$$Hv=nv.$$
We call $n$ the {\it weight} of the representation. Moreover, the
set $\{v, Yv,\ldots Y^nv\}$ is a basis of $V$, and we have
\begin{eqnarray*}
H(Y^iv)&=& (n-2i)Y^iv\;(i=0,1,\ldots, n),\\
X(Y^iv)&=& i(n-i+1)Y^{i-1}v\;(i=0,1,\ldots, n),\\
Y(Y^iv)&=& Y^{i+1}v\; (i=0,1,\ldots, n-1),\\
Y(Y^nv)&=&0.
\end{eqnarray*}

\bigskip
\noindent {\bf Remark 2.2.} The trivial representation
$V_0=\overline {\bf Q}_l$ of $\hbox{\frak sl}(2)$ is the irreducible
representation of weight $0$. Let $V_1=\overline {\bf Q}_l^2$ be the
standard representation of $\hbox{\frak sl}(2)$ on which
$\hbox{\frak sl}(2)$ acts as the multiplication of matrices on
column vectors. It is the irreducible representation of weight $1$,
and $f_0=\left(\begin{array}{c}1\\0\end{array}\right)$ is a highest
weight vector. Let $V_n={\rm Sym}^n(V_1)$ be the $n$-th symmetric
product of $V_1$. It is the irreducible representation of weight
$n$, and $f_0^n$ is a highest weight vector.

\bigskip
Let $V_n$ be the irreducible representation of $\hbox{\frak sl}(2)$
of weight $n$. Note that the eigenvalues $n,n-2, n-4, \ldots, -n$ of
$H$ form an unbroken arithmetic progression of integers with
difference $-2$, and each eigenvalue has multiplicity 1. Moreover,
the space ${\rm ker}(X)$ has dimension 1 and coincides with the
eigenspace of $H$ corresponding to the eigenvalue $n$. For any
integer $w$, let $V_n^w$ be the eigenspace of $H$ corresponding to
the eigenvalue $w$. We then have
$${\rm dim}(V_n^w)=\left\{\begin{array}{cl}
1&\hbox { if } w\equiv n \hbox { mod }2 \hbox { and } -n\leq w\leq n,\\
0&\hbox { otherwise.}\end{array}\right.$$ Moreover, we have
\begin{eqnarray*}
V_n\cap{\rm ker}(X)&=&V_n^n, \\
V_n^w\cap {\rm ker}(X)&=&\left\{\begin{array}{cl}
V_n^w&\hbox { if } w=n,\\
0&\hbox { otherwise.}\end{array}\right.
\end{eqnarray*}
In general, any finite dimensional representation $V$ of
$\hbox{\frak sl}(2)$ is a direct sum of irreducible representations.
Let
$$V=m_0V_0\oplus m_1V_1\oplus \cdots\oplus m_kV_k$$ be the isotypic
decomposition of $V$. For any integer $w$, let $V^w$ be the
eigenspace of $H$ corresponding to the eigenvalue $w$. If $w$ is
non-negative, then we have
$$V^w=m_wV_w^w\oplus m_{w+2}V_{w+2}^w\oplus\cdots$$
and
$${\rm dim}(V^w)=m_w+m_{w+2}+\cdots.$$
Moreover, we have
\begin{eqnarray*}
{\rm ker}(X)&=&(m_0V_0\cap {\rm ker}(X))\oplus (m_1V_1\cap{ \rm
ker}(X))\oplus \cdots \oplus (m_kV_k\cap{\rm ker}(X))\\
&=& m_0V_0^0\oplus m_1V_1^1\oplus\cdots\oplus m_kV_k^k.
\end{eqnarray*}
and hence
$${\rm ker}(X)\cap V^w= m_wV_w^w.$$ It follows that
$${\rm ker}X=({\rm ker}(X)\cap V^0)\oplus ({\rm ker}(X)\cap V^1)
\oplus\cdots\oplus ({\rm ker}(X)\cap V^k)$$ and
$${\rm dim}({\rm ker}(X)\cap V^w)= m_w\\
={\rm dim}(V^w)-{\rm dim}(V^{w+2}).$$ We summarize these results as
follows.

\bigskip
\noindent {\bf Proposition 2.3.} Let $V$ be a finite dimensional
$\overline {\bf Q}_l$-representation of $\hbox{\frak sl}(2)$. For
any integer $w$, let $V^w$ be the eigenspace of $H$ corresponding to
the eigenvalue $w$. Then we have
$${\rm ker}X=({\rm ker}(X)\cap V^0)\oplus ({\rm ker}(X)\cap V^1)
\oplus\cdots,$$ and for any non-negative $w$, we have
$${\rm dim}({\rm ker}(X)\cap V^w)
={\rm dim}(V^w)-{\rm dim}(V^{w+2}).$$

\bigskip
\bigskip
\centerline {\bf 3. The Local Factor at $0$}

\bigskip
\bigskip
In this section, we calculate the local factor
$${\rm det}(I-F_0t, ({\rm Sym}^k({\rm Kl}_{n+1}))^{I_0})$$ at $0$ of
the $L$-function of the $k$-th symmetric product of the Kloosterman
sheaf. Let $(V, N, F)$ be the triple defined in \S 1 corresponding
to the Kloosterman sheaf. Then the above local factor is simply
$${\rm det}\left(I-Ft, {\rm ker}(N: {\rm Sym}^k(V)\to {\rm
Sym}^k(V))\right).$$

\bigskip
Let $V_1=\overline {\bf Q}_l^2$ be the standard representation of
$\hbox{\frak sl}(2)$. Set
$$f_0=\left(\begin{array}{c}1\\0\end{array}\right),\;
f_1=\left(\begin{array}{c}0\\1\end{array}\right).$$ We have
\begin{eqnarray*}
&&H(f_0)=f_0,\; H(f_1)=-f_1,\\
&&X(f_0)=0,\; X(f_1)=f_0.
\end{eqnarray*}
Let $V_{n-1}={\rm Sym}^{n-1}(V_1)$, and set
$$e_i=\frac{1}{i!}f_0^{n-1-i}f_1^{i}\; (i=0,1,\ldots, n-1).$$ We have
$$H(e_0)=(n-1)e_0,\; H(e_1)=(n-3)e_1,\; \ldots, \;H(e_{n-1})=-(n-1)e_{n-1}$$ and
$$X(e_0)=0,\; X(e_1)=e_0,\; \ldots,\; X(e_{n-1})=e_{n-2}.$$
Comparing with Proposition 1.1, we can identify $V_{n-1}$ with $V$
coming from the triple $(V,F,N)$ defined by the Kloosterman sheaf
such that $N$ is identified with $X$, and the eigenspace of $F$ with
eigenvalue $q^i$ is identified with the eigenspace of $H$ with
eigenvalue $n-2i-1$.

Consider the $k$-th symmetric product ${\rm Sym}^k(V_{n-1})$. It has
a basis $$\{e_0^{i_0}e_1^{i_1}\cdots e_n^{i_{n-1}}|i_j\geq 0,
i_0+i_1+\cdots +i_{n-1}=k\}.$$ We have
$$H(e_0^{i_0}e_1^{i_1}\cdots e_{n-1}^{i_{n-1}})=((n-1)\cdot i_0+(n-3)\cdot i_1+\cdots+
(-(n-1))\cdot i_{n-1})e_0^{i_0}e_1^{i_1}\cdots  e_{n-1}^{i_{n-1}}.$$
So $e_0^{i_0}e_1^{i_1}\cdots e_{n-1}^{i_{n-1}}$ is an eigenvector of
$H$ with eigenvalue $(n-1)\cdot i_0+(n-3)\cdot i_1+\cdots+
(-(n-1))\cdot i_{n-1}$. It is also an eigenvector $F$ with
eigenvalue
$$q^{0\cdot i_0+1\cdot i_1+\cdots +(n-1)\cdot i_{n-1}}
=q^{\frac{1}{2}((n-1)k-((n-1)\cdot i_0+(n-3)\cdot i_1+\cdots+
(-(n-1))\cdot i_{n-1}))}.$$ Here we use the fact that
\begin{eqnarray*}
&&2(0\cdot i_0+1\cdot i_1+\cdots +(n-1)\cdot i_{n-1})+((n-1)\cdot
i_0+(n-3)\cdot i_1+\cdots+ (-(n-1))\cdot
i_{n-1})\\
&=&(n-1)(i_0+i_1+\cdots+ i_{n-1})\\
&=&k(n-1).
\end{eqnarray*}
This equality also shows that $$(n-1)\cdot i_0+(n-3)\cdot
i_1+\cdots+ (-(n-1))\cdot i_{n-1}\equiv k(n-1) \hbox { mod } 2.$$
For each non-negative integer $w$, let
$$D_k(w)=\{(i_0,\ldots, i_{n-1})| i_j\geq 0, \;
i_0+i_1+\cdots +i_{n-1}=k,\;(n-1)\cdot i_0+(n-3)\cdot i_1+\cdots+
(-(n-1))\cdot i_{n-1}=w\},$$ and let $d_k(w)$ be the number of
elements of $D_k(w)$. We have $d_k(w)=0$ if $w\not \equiv k(n-1)
\hbox { mod }2$ or if $w> k(n-1)$. Note that
$$\{e_0^{i_0}e_1^{i_1}\cdots e_{n-1}^{i_{n-1}}|(i_0,i_1,\ldots, i_{n-1})\in
D_k(w)\}$$ is a basis of the eigenspace $({\rm Sym}^k(V_{n-1}))^w$
of $H$ with eigenvalue $w$. By Proposition 2.3, we have
$${\rm ker}(X)=\bigoplus_{w=0}^{k(n-1)} {\rm ker}(X)\cap ({\rm
Sym}^k(V_{n-1}))^w$$ and
$${\rm dim}({\rm ker}(X)\cap ({\rm
Sym}^k(V_{n-1}))^w)=d_k(w)-d_k(w+2).$$ Now $({\rm
Sym}^k(V_{n-1}))^w$ is also the eigenspace of $F$ on ${\rm
Sym}^k(V)$ with eigenvalue $q^{\frac{k(n-1)-w}{2}}$. So we have
$${\rm det}\left(I-Ft, {\rm ker}(N: {\rm Sym}^k(V)\to {\rm
Sym}^k(V))\right)=\prod_{w=0}^{k(n-1)}
(1-q^{\frac{k(n-1)-w}{2}}t)^{d_k(w)-d_k(w+2)}.$$ As  $d_k(w)=0$ if
$w\not \equiv k(n-1) \hbox { mod }2$ or if $w>k (n-1)$, we have
$${\rm det}\left(I-Ft, {\rm ker}(N: {\rm Sym}^k(V)\to {\rm
Sym}^k(V))\right)=\prod_{u=0}^{[\frac{k(n-1)}{2}]}
(1-q^ut)^{d_k(k(n-1)-2u)-d_k(k(n-1)-2u+2)}.$$ Set
$$c_k(u)=d_k(k(n-1)-2u)$$ so that we have
$${\rm
det}\left(I-Ft, {\rm ker}(N: {\rm Sym}^k(V)\to {\rm
Sym}^k(V))\right)=\prod_{u=0}^{[\frac{k(n-1)}{2}]}
(1-q^ut)^{c_k(u)-c_k(u-1)}.$$ In the following, we find an
expression for $c_k(u)-c_k(u-1)$.

Note that $c_k(u)$ is the number of elements of the set
$$\{(i_0,\ldots, i_{n-1})| i_j\geq 0, \;
i_0+i_1+\cdots +i_{n-1}=k,\;0\cdot i_0+1\cdot i_1+\cdots +(n-1)\cdot
i_{n-1}=u\}.$$ Taking power series expansion, we get
$$\frac{1}{(1-y)(1-xy)\cdots (1-x^{n-1}y)}=\sum_{k=0}^\infty
\sum_{u=0}^\infty c_k(u) x^uy^k.$$ Since
$$(1-y)\frac{1}{(1-y)(1-xy)\cdots (1-x^{n-1}y)}=(1-x^{n}y)
\frac{1}{(1-xy)\cdots (1-x^{n}y)},$$ we have
$$(1-y)(\sum_{k,u}c_k(u)x^uy^k)=(1-x^{n}y)(\sum_{k,u}c_k(u)x^u(xy)^k),$$
that is,
$$\sum_{k,u}c_k(u)x^uy^k-\sum_{k,u}c_k(u)x^uy^{k+1}
=\sum_{k,u}c_k(u)x^{u+k}y^k-\sum_{k,u}c_k(u)x^{n+u+k}y^{k+1}.$$
Comparing the coefficients of $y^k$, we get
$$\sum_u c_k(u)x^u-\sum_u c_{k-1}(u)x^u=(\sum_uc_k(u)x^u)x^k-
(\sum_u c_{k-1}(u)x^u)x^{n+k-1},$$ that is,
$$\sum_u c_k(u)x^u=\frac{1-x^{n+k-1}}{1-x^k}\sum_u c_{k-1}(u)x^u.$$
Applying this expression repeatedly, we get
\begin{eqnarray*}
\sum_u c_k(u)x^u&=&\frac{1-x^{n+k-1}}{1-x^k}\sum_u c_{k-1}(u)x^u\\
&=&\frac{(1-x^{n+k-2})(1-x^{n+k-1})}{(1-x^{k-1})(1-x^k)}\sum_u c_{k-2}(u)x^u\\
&=&\cdots\\
&=&\frac{(1-x^{n})\cdots (1-x^{n+k-2})(1-x^{n+k-1})}{(1-x)\cdots
(1-x^{k-1})(1-x^k)}.
\end{eqnarray*}
Therefore
\begin{eqnarray*}
\sum_u(c_k(u)-c_k(u-1))x^u&=& \sum_u c_k(u)x^u-x\sum_u c_k(u)x^u\\
&=& (1-x)\sum_u c_k(u)x^u\\
&=&(1-x)\frac{(1-x^{n})\cdots
(1-x^{n+k-2})(1-x^{n+k-1})}{(1-x)\cdots (1-x^{k-1})(1-x^k)}\\
&=& \frac{(1-x^{n})\cdots (1-x^{n+k-2})(1-x^{n+k-1})}{(1-x^2)\cdots
(1-x^{k-1})(1-x^k)}.
\end{eqnarray*}
So $c_k(u)-c_k(u-1)$ is the coefficients of $x^u$ in the power
series expansion of $\frac{(1-x^{n})\cdots
(1-x^{n+k-2})(1-x^{n+k-1})}{(1-x^2)\cdots (1-x^{k-1})(1-x^k)}$. We
finally get the following, which is Theorem 0.2 in the Introduction.

\bigskip
\noindent {\bf Theorem 3.1.} We have
$${\rm det}(I-F_0t, ({\rm Sym}^k({\rm Kl}_{n}))^{I_0})=\prod_{u=0}^
{[\frac{k(n-1)}{2}]} (1-q^ut)^{m_k(u)},$$ where $m_k(u)$ is
determined by
$$\frac{(1-x^{n})\cdots (1-x^{n+k-2})(1-x^{n+k-1})}{(1-x^2)\cdots
(1-x^{k-1})(1-x^k)}=\sum_{u=0}^\infty m_k(u)x^u.$$ We have
$$m_k(u)=c_k(u)-c_k(u-1),$$ where $c_k(u)$ is the number of elements of the set
$$\{(i_0,\ldots, i_{n-1})| i_j\geq 0, \;
i_0+i_1+\cdots +i_{n-1}=k,\;0\cdot i_0+1\cdot i_1+\cdots +(n-1)\cdot
i_{n-1}=u\}.$$

\bigskip
\bigskip

\centerline {\bf 4. $L$-functions of $p$-Adic Symmetric Products}

\bigskip
\bigskip
Let $s$ be a $p$-adic integer. Define the $p$-adic symmetric product
$L$-function to be the $p$-adic limit
$$L_p(s,n, T) = \lim_{i\rightarrow \infty} L({\rm Sym}^{k_i}({\rm Kl}_{n}), T),$$
where $k_i$ is any sequence of increasing positive integers going to
infinity as complex numbers and approaching to $s$ as $p$-adic
integers. This $L$-function is really the $L$-function of some
infinite rank overconvergent nuclear $\sigma$-module. Confer [W2].
As a consequence, it is a $p$-adic meromorphic function. This
$L$-function plays the key role in the proof [W1] of Dwork's unit
root conjecture for the Kloosterman family. In this section, we
prove the following more precise results about $L_p(s,n, T)$.
\bigskip

\noindent {\bf Theorem 4.1.} Let $d_j$ be the coefficient of $x^j$
in the power series
$${1 \over (1-x^2)(1-x^3)\cdots (1-x^{n-1})}.$$
For each $p$-adic integer $s$, we can write
$$L_p(s,n,T) = A_p(s,n,T) \prod_{j=0}^{\infty} (1-q^jT)^{d_j},$$
where $A_p(s,n, T)$ is a $p$-adically entire function. In
particular, the $p$-adic series $L_p(s,n, T)$ is $p$-adically entire
and it has a zero at $T= q^{-j}$ with multiplicity at least $d_j$
for each non-negative integer $j$.

\bigskip
\noindent {\bf Proof}. Take a sequence $k_i$ of increasing positive
integers going to infinity as complex numbers and approaching to $s$
as $p$-adic integers. Since ${\rm Kl}_n$ is pure of weight $n-1$,
for each positive integer $k_i$, Grothendieck's formula for
$L$-functions implies that we can write
$$L(k_i,n,T):=L({\bf A}_{{\bf F}_q}^1, {\rm Sym}^{k_i}({\rm Kl}_{n}), T) = {P(k_i,n, T) \over
((1-q^{(n-1)k_i/2}T)(1-q^{((n-1)k_i+2)/2}T))^{e_i}},$$ where
$$P(k_i,n, T)={\rm det}(1-FT, H^1({\bf P}_{\bf
F}^1,j_\ast({\rm Sym}^k({\rm Kl}_n)))){\rm det}(1-F_0T,(({\rm
Sym}^k{\rm Kl}_n)_{\bar \eta})^{I_0}){\rm det}(1-F_\infty T,(({\rm
Sym}^k{\rm Kl}_n)_{\bar \eta})^{I_\infty}),
$$  and $e_i$ is the
multiplicity of the geometrically trivial representation in ${\rm
Sym}^{k_i}({\rm Kl}_{n})$. In fact, by Corollary 0.5, we know that
$e_i=0$ unless $p=2$, $k_i$ even and $n$ odd, in which case we have
$e_i=1$. Taking the limit, we deduce that
$$L_p(s,n, T) = \lim_{i\rightarrow \infty} P(k_i,n, T).$$
Fix a positive integer $r$. By the results in [W1] (Theorem 5.7 and
Lemma 5.10), the number of zeros and poles of the $L$-function
$L(k_i, n,T)$ as $k_i $ varies is uniformly bounded in the disk
$|T|_p<p^r$. In particular, the number of zeros of the polynomial
$P(k_i, n,T)$ (the numerator of $L(k_i, n, T)$) as $k_i $ varies is
uniformly bounded in the disk $|T|_p<p^r$. Under the condition
$k\geq n$, we have
$$\frac{(1-x^{n})\cdots (1-x^{n+k-2})(1-x^{n+k-1})}{(1-x^2)\cdots
(1-x^{k-1})(1-x^k)}=\frac{(1-x^{k+1})\cdots
(1-x^{n+k-2})(1-x^{n+k-1})}{(1-x^2)\cdots (1-x^{n-2})(1-x^{n-1})}.$$
It follows that $m_k(j)=d_j$ for $j\leq k$, where $m_k(j)$ is
defined in Theorem 3.1. So we have $m_{k_i}(j) = d_j$  for all
$1\leq j\leq r$ provided that $k_i\geq {\rm max}(r,n)$. Then by
Theorem 3.1, we can write
$$P(k_i,n, T) = B_r(k_i,n, T) \prod_{j=0}^r
(1-q^jT)^{d_j},$$ where $B_r(k_i, n,T)\in 1+T{\bf Z}[T]$ is a
polynomial in $T$. Furthermore, the number of the zeros of
$B_r(k_i,n, T)$ in the disk $|T|_p <p^r$ is uniformly bounded as
$k_i$ varies. This implies that the limit
$$C_r(s,n, T): =\lim_{i\rightarrow \infty} B_r(k_i,n, T) =
{L_p(s,n, T) \over \prod_{j=0}^r (1-q^jT)^{d_j}}$$exists and is
$p$-adically analytic in the disk $|T|_p <p^r$. In particular,
$L_p(s,n, T)$ is $p$-adically analytic in the disk $|T|_p <p^r$ and
has a zero at $T=q^{-j}$ with multiplicity at least $d_j$ for $0\leq
j\leq r$. As we can take $r$ to be an arbitrarily large integer, we
deduce that
$$A_p(s,n,T):= \lim_{r\rightarrow \infty} C_r(s,n, T)=
{L_p(s,n, T) \over \prod_{j=0}^{\infty} (1-q^jT)^{d_j}}$$ is
$p$-adically entire. The theorem is proved.

\bigskip
Note that Theorem 0.2 shows that for $(n,p)=1$, the limit of the
local factors at infinity disappears and hence has no contribution
to the zeros of $L$-functions of $p$-adic infinite symmetric
products. This together with the above proof implies that
$$A_p(s,n,T) = \lim_{i\rightarrow \infty}
{\rm det}(1-FT, H^1({\bf P}_{\bf F}^1,j_\ast({\rm Sym}^{k_i}({\rm
Kl}_n)))),
$$
that is, $A_p(s,n,T)$ is the $p$-adic limit of the non-trivial
factor of $L({\bf A}_{{\bf F}_q}^1-\{0\},{\rm Sym}^{k_i}({\rm
Kl}_{n}), T)$ as $k_i$ approaches to $s$. It is a $p$-adic entire
function. Its zeros are called non-trivial zeros of $L_p(s,n,T)$.
Some partial results on the distribution of the zeros of $L_p(s,n,
T)$ were obtained in [W2].

\bigskip
\noindent {\bf Remark.} The same proof shows that the entireness
property for $L_p(s,n,T)$ can be extended to any $p$-adic $s$-th
symmetric product $L$-function of a lisse pure positive weight
$l$-adic sheaf whose $p$-adic unit part has rank one and is a
$p$-adic $1$-unit. The Kloosterman sheaf is just the first such
example. The ordinary family of Calabi-Yau hypersurfaces is another
important example, generalizing the ordinary family of elliptic
curves which has been well studied in connection with the theory of
$p$-adic modular forms.

\bigskip
\bigskip
\centerline {\bf 5. Variation of the non-trivial factor}

\bigskip
\bigskip
In this section, we derive some consequences for the non-trivial
factor
$$K_q(k,n,T):={\rm det}(1-FT, H^1({\bf P}_{\bf F}^1,j_\ast({\rm
Sym}^k({\rm Kl}_n))))\in 1+T{\bf Z}[T].$$ This is a polynomial with
integer coefficients, pure of weight $k(n-1)+1$. Its degree can be
computed explicitly by the degree formula for $L(k,n,T)$ (Theorem
0.1 in [FW]) and the degree formulas for the trivial factors of
$L(k,n,T)$ as implicit in Theorem 0.2, Theorem 0.3 and Corollary
0.5.

In the simplest case $n=2$ and $q=p$,  the polynomial $K_p(k,n,T)$
is the Kloosterman analogue of the $p$-th Hecke polynomial acting on
weight $k+2$ modular forms. It would be interesting to understand
how the polynomial $K_p(k,n,T)$ varies as $p$ varies while $k$ is
fixed or as $k$ varies while $p$ is fixed.

For fixed $k$ and $n$, the polynomial $K_p(k,n,T)$ should be the
$p$-th Euler factor of a motive $M_{k,n}$ over ${\bf Q}$. It would
be interesting to construct explicitly this motive (its underlying
scheme) or its corresponding compatible system of Galois
representation or its automorphic interpretation. In the special
case when $n=2$ and $k=5,6$, the polynomial $K_p(k,n,T)$ has degree
$2$ and is conjectured by Choi-Evans-Stark ([CE]) to be the Euler
factor at $p$ of an explicit modular form of weight $k+2$.

Just like the case for $L(k,n,T)$, we are interested in how the
polynomial $K_p(k,n,T)$ varies as $k$ varies $p$-adically. This
question was studied by Gouvea-Mazur for $p$-th Hecke polynomials in
connection with $p$-adic variation of modular forms. The first
simple result is a $p$-adic continuity result.

\bigskip
\noindent {\bf Proposition 5.1}. Let $k_1, k_2$ and $k_3$ be
positive integers such that $k_1 = k_2 +p^mk_3$ with $k_1$ not
divisible by $p$. Then we have the congruence
$$K_p(k_1,n,T) \equiv K_p(k_2,n,T) ({\rm mod}~ p^{\min (m,
k_2/2)}).$$

\bigskip
\noindent{\bf Proof}: Let $q=p$. The Frobenius eigenvalues of the
Kloosterman sheaf at each closed point are all divisible by $p$
except for exactly one eigenvalue which is a $p$-adic $1$-unit. From
this and the Euler product definition of the $L$-function
$L(k,n,T)$, we deduce the slightly stronger congruence:
$$L(k_1,n,T) \equiv L(k_2,n,T) ({\rm mod}~ p^{\min (m,
k_2)}).$$ To prove the proposition, it remains to check that the
same congruence in the proposition holds for the trivial factors.
This follows from the explicit results stated in Theorem 0.2,
Theorem 0.3 and Corollary 0.5.

\bigskip
Let $s$ be a $p$-adic integer. Choose a sequence of positive
integers $k_i$ going to infinity as complex numbers and
approaching $s$ as $p$-adic integers. The above congruence for
$K_p(k,n,T)$ implies that the limit
$$A_p(s,n,T):= \lim_{i\rightarrow \infty} K_p(k_i,n,T)$$
exists and it is exactly the non-trivial factor $A_p(s,n,T)$ in
Theorem 4.1. It follows that $A_p(s,n,T)$ is a $p$-adic entire
function. It would be interesting to determine the $p$-adic Newton
polygon of the entire function $A_p(s,n,T)$. This would give exact
information on the distribution of the zeros of $A_p(s,n,T)$.

The rigid analytic curve in the $(s,T)$ plane defined by the
equation $A_p(s,n,T)=0$ is the Kloosterman sum analogue of the
eigencurve in the theory of $p$-adic modular forms studied by
Coleman-Mazur [CM]. It would be interesting to study the properties
of the rigid analytic curve $A_p(s,n,T)=0$ and its relation to
$p$-adic automorphic forms.

\bigskip
\bigskip
\centerline {\bf 6. Appendix}

\bigskip
\bigskip
In this section, we prove the following proposition, which is Lemma
0.4 in the Introduction.

\bigskip
\noindent {\bf Proposition 6.1.} Let $\hbox {\frak g}$ be one of the
following Lie algebras
$$\hbox{\frak sl}(n), \hbox {\frak sp}(n),
\hbox{\frak so}(n), \hbox{\frak g}_2$$ and let $V$ be the standard
representation of $G$. Then in the case where $\hbox{\frak g}=\hbox
{\frak sl}(n)$ or $\hbox {\frak sp}(n)$, the representation ${\rm
Sym}^k(V)$ is irreducible, and in the case where $\hbox{\frak
g}=\hbox{\frak so}(n)$ or $\hbox{\frak g}_2$, the representation
${\rm Sym}^k(V)$ contains exactly one copy of the trivial
representation if $k$ is even, and contains no trivial
representation if $k$ is odd.

\bigskip
\noindent {\bf Proof.} First recall the dimension formula for
irreducible representations of simple Lie algebras. Let $\hbox
{\frak g}$ be a simple Lie algebra (over $\overline {\bf Q}_l$).
Choose a Cartan subalgebra $\hbox {\frak h}$ of $\hbox {\frak g}$,
and let $R$ be the set of roots. We have the Cartan decomposition
$$\hbox{\frak g}=\hbox{\frak h}\bigoplus(\bigoplus_{\alpha\in R} \hbox
{\frak g}_\alpha).$$ For each $\alpha\in R$, let $H_\alpha$ be the
unique element in $[\hbox {\frak g}_\alpha, \hbox {\frak
g}_{-\alpha}]$ such that $\alpha(H_\alpha)=2$. The weight lattice
$\Lambda_W$ is the lattice in $\hbox{\frak h}^\ast$ generated by
those linear functionals $\beta$ with the property
$\beta(H_\alpha)\in {\bf Z}$ for all $\alpha\in R$. Fix an ordering
of $R$. Let $R^+$ be the set of positive roots, and let ${\cal W}$
be the Weyl chamber. Set
$$\rho=\frac{1}{2}\sum_{\alpha\in R^+}\alpha.$$ For any $\lambda\in
\Lambda_W\cap {\cal W}$, the dimension of the irreducible
representation $\Gamma_\lambda$ with highest weight $\lambda$ is
given by
$${\rm dim}(\Gamma_\lambda)=\prod_{\alpha\in R^+} \frac{\langle \lambda
+\rho,\alpha\rangle}{\langle \rho,\alpha\rangle}=\prod_{\alpha\in
R^+} \frac{( \lambda +\rho,\alpha)}{(\rho,\alpha)},$$ where
$(\;,\;)$ is the Killing form on $\hbox{\frak h}^\ast$, and
$$\langle \beta,\alpha\rangle=\beta(H_\alpha)=\frac{2(\beta,\alpha)}
{(\alpha,\alpha)}$$ for any $\beta\in \hbox{\frak h}^\ast$ and
$\alpha\in R$.

For each pair $1\leq i, j\leq n$, let $E_{ij}$ be the $(n\times
n)$-matrix whose only nonzero entry is on the $i$-th row and $j$-th
column, and this nonzero entry is $1$. For each $1\leq i\leq n$, let
$L_i$ be the linear functional on the space of diagonal matrices
with the property
$$L_i(E_{jj})=\left\{\begin{array}{cl}
1& \hbox { if } i=j,\\
0&\hbox { if } i\not= j.\end{array}\right.$$

Consider the Lie algebra $\hbox{\frak sl}(n)$ of traceless $(n\times
n)$-matrices. Let $\hbox{\frak h}$ be the space of diagonal matrices
in $\hbox{\frak sl}(n)$. It is a Cartan subalgebra of $\hbox{\frak
sl}(n)$. The set of roots of $\hbox{\frak sl}(n)$ are
$$R=\{L_i-L_j|i\not = j\}$$ and $$H_{L_i-L_j}=E_{ii}-E_{jj}\; (i\not=j).$$
Choose an ordering of roots so that
$$R^+=\{L_i-L_j|i<j\}$$ is the set of the positive roots.
We have
$$\rho=\sum_{i=1}^n (n-i)L_i.$$
(To deduce this formula, we use the fact that $L_1+\cdots +L_n=0$
for $\hbox{\frak sl}(n)$.) By the dimension formula, for any
$$\lambda=\lambda_1L_1+\cdots +\lambda_nL_n$$ lying in the
intersection of the weight lattice and the Weyl chamber, the
dimension of the irreducible representation $\Gamma_{\lambda}$ of
$\hbox {\frak sl}(n)$ with highest weight $\lambda$ is
\begin{eqnarray*}
{\rm dim} (\Gamma_{\lambda})&=& \prod_{\alpha \in R^+}\frac{\langle
\lambda+\rho,
\alpha\rangle}{\langle \rho,\alpha\rangle}\\
&=&\prod_{i<j} \frac{(\lambda+\rho)(E_{ii}-E_{jj})}{\rho(E_{ii}-E_{jj})}\\
&=&\prod_{i< j}\frac{(\sum_i
(\lambda_i+(n-i)L_i))(E_{ii}-E_{jj})}{(\sum_i
(n-i)L_i)(E_{ii}-E_{jj})}\\
&=&\prod_{i<j}\frac{\lambda_i-\lambda_j+j-i}{j-i}
\end{eqnarray*}
In the case where $\lambda=kL_1$, we have
$$\lambda_i=\left\{\begin{array}{cl} k &\hbox { if } i=1,\\
0 &\hbox { if } i\geq 2.\end{array}\right.$$ So we have
$${\rm dim}(\Gamma_{kL_1})=\prod_{1<j}\frac{k+j-1}{j-1}=\left(\begin{array}{c}
k+n-1\\ n-1\end{array} \right).$$ Note that the dimension of
$\Gamma_{kL_1}$ is exactly the dimension of ${\rm Sym}^k(V)$. Since
the weights of the standard representation $V$ are $L_1,\ldots,
L_n$, the representation ${\rm Sym}^k (V)$ has a highest weight
$kL_1$. So we must have
$${\rm Sym}^k(V)=\Gamma_{kL_1}.$$ In particular, ${\rm Sym}^k(V)$ is
irreducible.

Now suppose $n=2m$ is an even number and consider the Lie algebra
$\hbox{\frak sp}(n)$ of matrices of the form
$$\left(\begin{array}{cc}
A&B\\
C&D \end{array}\right),$$ where $A,B,C,D$ are $(m\times
m)$-matrices, $B$ and $C$ are symmetric and  $A^t+D=0$. Let
$\hbox{\frak h}$ be the space of diagonal matrices in $\hbox{\frak
sp}(n)$. It is a Cartan subalgebra of $\hbox{\frak sp}(n)$. The set
of roots of $\hbox{\frak sp}(n)$ are
$$R=\{\pm L_i\pm L_j|1\leq i,j\leq m\}-\{0\}$$ and
\begin{eqnarray*}
H_{L_i-L_j}&=&(E_{ii}-E_{m+i,m+i})-(E_{jj}-E_{m+j,m+j})\; (i\not=j),\\
H_{L_i+L_j}&=&(E_{ii}-E_{m+i,m+i})+(E_{jj}-E_{m+j,m+j})\; (i\not=j),\\
H_{2L_i}&=&E_{ii}-E_{m+i,m+i}.
\end{eqnarray*}
Choose an ordering of roots so that
$$R^+=\{L_i-L_j|i<j\}\cup \{L_i+L_j|i\leq j\}$$ is the set of the positive roots.
We have
$$\rho=\sum_{i=1}^m (m+1-i)L_i.$$
By the dimension formula, for any $$\lambda=\lambda_1L_1+\cdots
+\lambda_mL_m$$ lying in the intersection of the weight lattice and
the Weyl chamber, the dimension of the irreducible representation
$\Gamma_{\lambda}$ of $\hbox {\frak sp}(n)$ with highest weight
$\lambda$ is
\begin{eqnarray*}
{\rm dim}(\Gamma_{\lambda})&=& \prod_{\alpha \in R^+}\frac{\langle
\lambda+\rho,
\alpha\rangle}{\langle \rho,\alpha\rangle}\\
&=&\prod_{i<j}
\frac{(\lambda+\rho)((E_{ii}-E_{m+i,m+i})-(E_{jj}-E_{m+j,m+j}))}
{\rho((E_{ii}-E_{m+i,m+i})-(E_{jj}-E_{m+j,m+j}))} \\
&& \prod_{i<j}
\frac{(\lambda+\rho)((E_{ii}-E_{m+i,m+i})+(E_{jj}-E_{m+j,m+j}))}
{\rho((E_{ii}-E_{m+i,m+i})+(E_{jj}-E_{m+j,m+j}))} \\
&& \prod_i
\frac{(\lambda+\rho)(E_{ii}-E_{m+i,m+i})}{\rho(E_{ii}-E_{m+i,m+i})}\\
&=& \prod_{i<j} \frac{\lambda_i-\lambda_j+j-i} {j-i} \prod_{i<j}
\frac{\lambda_i+\lambda_j+2m+2-i-j} {2m+2-i-j} \prod_i
\frac{\lambda_i+m+1-i}{m+1-i}.
\end{eqnarray*}
In the case where $\lambda=kL_1$, we have
$$\lambda_i=\left\{\begin{array}{cl} k &\hbox { if } i=1,\\
0 &\hbox { if } i\geq 2.\end{array}\right.$$ So we have
\begin{eqnarray*}
{\rm dim} (\Gamma_{kL_1})&=&\left(\prod_{1<j\leq
m}\frac{k+j-1}{j-1}\right)\left(\prod_{1<j\leq
m}\frac {k+2m+2-1-j}{2m+2-1-j}\right)\left( \frac{k+m+1-1}{m+1-1}\right)\\
&=&\frac{(k+1)(k+2)\cdots(k+m-1)}{1\cdot 2\cdots
(m-1)}\frac{(k+m+1)\cdots (k+2m-1)}{(m+1)\cdots (2m-1)}
\frac{k+m}{m}\\
&=& \left(\begin{array}{c} k+2m-1\\ 2m-1\end{array} \right).
\end{eqnarray*}
Note that the dimension of $\Gamma_{kL_1}$ is exactly the dimension
of ${\rm Sym}^k(V)$. Since the weights of the standard
representation $V$ are $L_1,\ldots, L_n$, the representation ${\rm
Sym}^k (V)$ has a highest weight $kL_1$. So we must have
$${\rm Sym}^k(V)=\Gamma_{kL_1}.$$ In particular, ${\rm Sym}^k(V)$
is irreducible.

Now consider the cases where $\hbox{\frak g}=\hbox {\frak so}(n)$ or
$\hbox {\frak g}_2$. In these cases, there is a symmetric
non-degenerate $\hbox{\frak g}$-invariant bilinear form $Q(\;,\;)$
on $V$. Consider the contraction map
\begin{eqnarray*}
{\rm Sym}^k(V)&\to &{\rm Sym}^{k-2}(V),\\
v_1\cdots v_k&\mapsto& \sum_{i<j}Q(v_i,v_j)v_1\cdots \hat
v_i\cdots\hat v_j\cdots v_k.
\end{eqnarray*}
It is an epimorphism of representations of $\hbox {\frak g}$. We
will show the kernel of the contraction map is irreducible.

First consider the case where $n=2m$ is even, and the Lie algebra is
$\hbox{\frak so}(n)$ of matrices of the form
$$\left(
\begin{array}{cc}
A&B\\
C&D
\end{array}
\right),$$ where $A,B,C,D$ are $(m\times m)$-matrices, $B$ and $C$
are skew-symmetric and $A^t+D=0$. Let $\hbox{\frak h}$ be the space
of diagonal matrices in $\hbox{\frak so}(n)$. It is a Cartan
subalgebra of $\hbox{\frak so}(n)$. The set of roots of $\hbox{\frak
so}(n)$ are
$$R=\{\pm L_i\pm L_j|1\leq i,j\leq m,\; i\not =j\}$$ and
\begin{eqnarray*}
H_{L_i-L_j}&=&(E_{ii}-E_{m+i,m+i})-(E_{jj}-E_{m+j,m+j})\; (i\not=j),\\
H_{L_i+L_j}&=&(E_{ii}-E_{m+i,m+i})+(E_{jj}-E_{m+j,m+j})\; (i\not=j).
\end{eqnarray*}
Choose an ordering of roots so that
$$R^+=\{L_i-L_j|i<j\}\cup \{L_i+L_j|i< j\}$$ is the set of the positive roots.
We have
$$\rho=\sum_{i=1}^m (m-i)L_i.$$
By the dimension formula, for any $$\lambda=\lambda_1L_1+\cdots
+\lambda_mL_m$$ lying in the intersection of the weight lattice and
the Weyl chamber, the dimension of the irreducible representation
$\Gamma_{\lambda}$ of $\hbox {\frak so}(n)$ with highest weight
$\lambda$ is
\begin{eqnarray*}
{\rm dim}(\Gamma_{\lambda})&=& \prod_{\alpha \in R^+}\frac{\langle
\lambda+\rho,
\alpha\rangle}{\langle \rho,\alpha\rangle}\\
&=&\prod_{i<j}
\frac{(\lambda+\rho)((E_{ii}-E_{m+i,m+i})-(E_{jj}-E_{m+j,m+j}))}
{\rho((E_{ii}-E_{m+i,m+i})-(E_{jj}-E_{m+j,m+j}))}\\
&&\prod_{i<j}
\frac{(\lambda+\rho)((E_{ii}-E_{m+i,m+i})+(E_{jj}-E_{m+j,m+j}))}
{\rho((E_{ii}-E_{m+i,m+i})+(E_{jj}-E_{m+j,m+j}))} \\
&=& \prod_{i<j} \frac{\lambda_i-\lambda_j+j-i} {j-i} \prod_{i<j}
\frac{\lambda_i+\lambda_j+2m-i-j} {2m-i-j}.
\end{eqnarray*}
In the case where $\lambda=kL_1$, we have
$$\lambda_i=\left\{\begin{array}{cl} k &\hbox { if } i=1,\\
0 &\hbox { if } i\geq 2.\end{array}\right.$$ So we have
\begin{eqnarray*}
{\rm dim}(\Gamma_{kL_1})&=&\prod_{1<j\leq
m}\frac{k+j-1}{j-1}\prod_{1<j\leq
m}\frac {k+2m-1-j}{2m-1-j} \\
&=&\frac{\frac{(k+m-1)!}{k!}}{(m-1)!}\frac{\frac{(k+2m-3)!}{(k+m-2)!}}
{\frac{(2m-3)!}{(m-2)!}}\\
&=& \frac{(k+m-1)(k+2m-3)!}{(m-1)(2m-3)!k!}\\
&=&\left(\begin{array}{c} k+2m-1\\ k\end{array}
\right)-\left(\begin{array}{c} k+2m-3\\ k-2\end{array} \right)\\
&=& {\rm dim}({\rm Sym}^k(V))-{\rm dim}({\rm Sym}^{k-2}(V)).
\end{eqnarray*}
Since the contraction map ${\rm Sym}^k(V)\to {\rm Sym}^{k-2}(V)$ is
surjective, and its kernel has a highest weight $kL_1$, it follows
that $\Gamma_{kL_1}$ coincides with the kernel of the contraction
map. So we must have
$${\rm Sym}^k(V)=\Gamma_{kL_1}\oplus {\rm Sym}^{k-2}(V).$$
Using this expression repeatedly, we get
$${\rm Sym}^k(V)=\bigoplus_{i=0}^{[\frac{k}{2}]}\Gamma_{(k-2i)L_1}.$$
In particular, when $k$ is even, ${\rm Sym}^k(V)$ contains one copy
of the trivial representation, and when $k$ is odd, it contains no
trivial representation.

Next consider the case where $n=2m+1$ is odd, and the Lie algebra is
$\hbox{\frak so}(n)$ of matrices of the form
$$\left(\begin{array}{ccc}
A&B&E\\
C&D&F\\
G&H&0\end{array} \right),$$ where $A,B,C,D$ are $(m\times
m)$-matrices, $E$ and $F$ are $(m\times 1)$-matrices, $G$ and $H$
are $(1\times m)$-matrices, $B$ and $C$ are skew-symmetric,
$A^t+D=0$, $E^t+H=0$, and $F^t+G=0$. Let $\hbox{\frak h}$ be the
space of diagonal matrices in $\hbox{\frak so}(n)$. It is a Cartan
subalgebra of $\hbox{\frak so}(n)$. The set of roots of $\hbox{\frak
so}(n)$ are
$$R=\{\pm L_i\pm L_j|1\leq i,j\leq m,\; i\not =j\}\cup
\{\pm L_i|1\leq i\leq m\}$$ and
\begin{eqnarray*}
H_{L_i-L_j}&=&(E_{ii}-E_{m+i,m+i})-(E_{jj}-E_{m+j,m+j})\; (i\not=j),\\
H_{L_i+L_j}&=&(E_{ii}-E_{m+i,m+i})+(E_{jj}-E_{m+j,m+j})\; (i\not=j),\\
H_{L_i}&=&2(E_{ii}-E_{m+i,m+i}).
\end{eqnarray*}
Choose an ordering of roots so that
$$R^+=\{L_i-L_j|i<j\}\cup \{L_i+L_j|i< j\}\cup \{L_i\}$$
is the set of the positive roots. We have
$$\rho=\sum_{i=1}^m (m+\frac{1}{2}-i)L_i.$$
By the dimension formula, for any $$\lambda=\lambda_1L_1+\cdots
+\lambda_mL_m$$ lying in the intersection of the weight lattice and
the Weyl chamber, the dimension of the irreducible representation
$\Gamma_{\lambda}$ of $\hbox {\frak so}(n)$ with highest weight
$\lambda$ is
\begin{eqnarray*}
{\rm dim}(\Gamma_{\lambda})&=& \prod_{\alpha \in R^+}\frac{\langle
\lambda+\rho,
\alpha\rangle}{\langle \rho,\alpha\rangle}\\
&=&\prod_{i<j}
\frac{(\lambda+\rho)((E_{ii}-E_{m+i,m+i})-(E_{jj}-E_{m+j,m+j}))}
{\rho((E_{ii}-E_{m+i,m+i})-(E_{jj}-E_{m+j,m+j}))}
\\&& \prod_{i<j}
\frac{(\lambda+\rho)((E_{ii}-E_{m+i,m+i})+(E_{jj}-E_{m+j,m+j}))}
{\rho((E_{ii}-E_{m+i,m+i})+(E_{jj}-E_{m+j,m+j}))}\\
&& \prod_i \frac{(\lambda+\rho)(2(E_{ii}-E_{m+i,m+i}))}
{\rho(2(E_{ii}-E_{m+i,m+i}))}
\\
&=& \prod_{i<j} \frac{\lambda_i-\lambda_j+j-i} {j-i} \prod_{i<j}
\frac{\lambda_i+\lambda_j+2m+1-i-j} {2m+1-i-j}\prod_i
\frac{\lambda_i+m+\frac{1}{2}-i}{m+\frac{1}{2}-i} .
\end{eqnarray*}
In the case where $\lambda=kL_1$, we have
$$\lambda_i=\left\{\begin{array}{cl} k &\hbox { if } i=1,\\
0 &\hbox { if } i\geq 2.\end{array}\right.$$ So we have
\begin{eqnarray*}
{\rm dim} (\Gamma_{kL_1})&=&\left(\prod_{1<j\leq
m}\frac{k+j-1}{j-1}\right)\left(\prod_{1<j\leq
m}\frac {k+2m+1-1-j}{2m+1-1-j}\right)\left(\frac{k+m+\frac{1}{2}-1}{m+\frac{1}{2}-1}\right) \\
&=&\frac{\frac{(k+m-1)!}{k!}}{(m-1)!}\frac{\frac{(k+2m-2)!}{(k+m-1)!}}
{\frac{(2m-2)!}{(m-1)!}}\frac{2k+2m-1}{2m-1}\\
&=& \frac{(2k+2m-1)(k+2m-2)!}{(2m-1)!k!}\\
&=&\left(\begin{array}{c} k+2m\\ k\end{array}
\right)-\left(\begin{array}{c} k+2m-2\\ k-2\end{array} \right)\\
&=&{\rm dim}({\rm Sym}^k(V))-{\rm dim}({\rm Sym}^{k-2}(V)).
\end{eqnarray*}
Since the contraction map ${\rm Sym}^k(V)\to {\rm Sym}^{k-2}(V)$ is
surjective, and its kernel has a highest weight $kL_1$, it follows
that $\Gamma_{kL_1}$ coincides with the kernel of the contraction
map. So we must have
$${\rm Sym}^k(V)=\Gamma_{kL_1}\oplus {\rm Sym}^{k-2}(V).$$
Using this expression repeatedly, we get
$${\rm Sym}^k(V)=\bigoplus_{i=0}^{[\frac{k}{2}]}\Gamma_{(k-2i)L_1}.$$
In particular, when $k$ is even, ${\rm Sym}^k(V)$ contains one copy
of the trivial representation, and when $k$ is odd, it contains no
trivial representation.

Finally let $n=7$ and consider the Lie algebra $\hbox{\frak g}_2$.
The following points on the real plane form the root system $R$ of
$\hbox {\frak g}_2$:
\begin{eqnarray*}
&&\alpha_1=(1,0)\\
&&\alpha_2=(\frac{3}{2},\frac{\sqrt {3}}{2}),\\
&&\alpha_3=(\frac{1}{2},\frac{\sqrt {3}}{2}),\\
&&\alpha_4=(0,\sqrt {3}),\\
&&\alpha_5=(-\frac{1}{2},\frac{\sqrt {3}}{2}),\\
&&\alpha_6=(-\frac{3}{2},\frac{\sqrt 3}{2}),\\
&&\beta_1=-\alpha_1,\;
\beta_2=-\alpha_2,\;\beta_3=-\alpha_3,\;\beta_4=-\alpha_4,\;\beta_5
=-\alpha_5,\;\beta_6=-\alpha_6.
\end{eqnarray*}
Moreover, the Killing form induces the canonical inner product on
the real plane spanned by the roots. Choose an order on $R$ so that
$\alpha_i$ $(i=1,\ldots, 6)$ are the positive roots. The Weyl
chamber ${\cal W}$ is the positive cone generated by $\alpha_3$ and
$\alpha_4$, and the weight lattice $\Lambda_W$ is the lattice
generated by $\alpha_1$ and $\alpha_6$. Any element in
$\Lambda_W\cap {\cal W}$ is of the form
$$\lambda=a\alpha_3+b\alpha_4=\left(\frac{1}{2}a, \frac{\sqrt 3}{2}a+\sqrt
3 b\right),$$ where $a$ and $b$ are non-negative integers. We have
$$\rho=\frac{1}{2}\sum_{i=1}^6 \alpha_i=\left(\frac{1}{2},
\frac{3\sqrt 3}{2}\right),$$ and
$$\begin{array}{cclccl}
(\lambda+\rho,\alpha_1)&=&\frac{1}{2}(a+1),&
(\rho,\alpha_1)&=&\frac{1}{2},\\
(\lambda+\rho,\alpha_2)&=&\frac{3}{2}(a+b+2),&
(\rho,\alpha_2)&=&3,\\
(\lambda+\rho,\alpha_3)&=&\frac{1}{2}(2a+3b+5),&
(\rho,\alpha_3)&=&\frac{5}{2},\\
(\lambda+\rho,\alpha_4)&=&\frac{3}{2}(a+2b+3),&
(\rho,\alpha_4)&=&\frac{9}{2},\\
(\lambda+\rho,\alpha_5)&=&\frac{1}{2}(a+3b+4),&
(\rho,\alpha_5)&=&2,\\
(\lambda+\rho,\alpha_6)&=&\frac{3}{2}(b+1),&
(\rho,\alpha_6)&=&\frac{3}{2}.
\end{array}$$
By the dimension formula, the dimension of the irreducible
representation $\Gamma_{\lambda}$ with highest weight
$\lambda=a\alpha_3+b\alpha_4$ is
\begin{eqnarray*}
&&{\rm dim}(\Gamma_{\lambda})\\&=&\prod_{\alpha\in
R^+}\frac{(\lambda+\rho, \alpha)}{(\rho,\alpha)}\\
&=&\frac{\frac{1}{2}(a+1)\cdot
\frac{3}{2}(a+b+2)\cdot\frac{1}{2}(2a+3b+5)\cdot
\frac{3}{2}(a+2b+3)\cdot \frac{1}{2}(a+3b+4)\cdot
\frac{3}{2}(b+1)}{\frac{1}{2}\cdot 3 \cdot \frac{5}{2}\cdot
\frac{9}{2}\cdot 2\cdot \frac{3}{2}}\\
&=&\frac{(a+1)(a+b+2)(2a+3b+5)(a+2b+3)(a+3b+4)(b+1)}{120}.
\end{eqnarray*}
In particular, the dimension of the irreducible representation
$\Gamma_{\alpha_3}$ is
$${\rm dim} (\Gamma_{\alpha_3})=\frac{2\cdot 3\cdot 7\cdot 4\cdot
5\cdot 1 }{120}=7.$$ So $\Gamma_{\alpha_3}$ is the standard
representation $V$. The dimension of the irreducible representation
$\Gamma_{k\alpha_3}$ is
\begin{eqnarray*}
{\rm dim}(\Gamma_{k\alpha_3})&=&
\frac{(k+1)(k+2)(2k+5)(k+3)(k+4)}{120}
\\
&=&\left(\begin{array}{c} k+6\\6
\end{array}\right)-\left(\begin{array}{c} k+4\\6
\end{array}\right)\\
&=&{\rm dim}({\rm Sym}^k(V))-{\rm dim}({\rm Sym}^{k-2}(V)).
\end{eqnarray*}
Since the contraction map ${\rm Sym}^k(V)\to {\rm Sym}^{k-2}(V)$ is
surjective, and its kernel has a highest weight $k\alpha_3$, the
representation $\Gamma_{k\alpha_3}$ coincides with the kernel of the
contraction map. So we must have
$${\rm Sym}^k(V)=\Gamma_{k\alpha_3}\oplus {\rm Sym}^{k-2}(V).$$
Using this expression repeatedly, we get
$${\rm Sym}^k(V)=\bigoplus_{i=0}^{[\frac{k}{2}]}\Gamma_{(k-2i)\alpha_3}.$$
In particular, when $k$ is even, ${\rm Sym}^k(V)$ contains one copy
of the trivial representation, and when $k$ is odd, it contains no
trivial representation. This finishes the proof of the proposition.

\bigskip
\bigskip
\noindent {\bf References}

\bigskip
\bigskip

\noindent [CE] H. T. Choi and R. Evans, {\it Congruences for sums of
powers of Kloosterman sums}, Int. J. Number Theory, in press (2006).
\bigskip

\noindent [CM] R. Coleman and B. Mazur, {\it The Eigencurve}, in
Galois Representations in Arithmetic Algebraic Geometry (Durham,
1996), 1--113, London Math. Soc. Lecture Note Ser., 254, Cambridge
Univ. Press, Cambridge, 1998.

\bigskip

\noindent [D1] P. Deligne, {\it Applications de la Formule des
Traces aux Sommes Trigonom\'etriques}, in Cohomologie \'Etale (SGA
$4\frac{1}{2}$), 168-232, Lecture Notes in Math. 569,
Springer-Verlag 1977.

\bigskip
\noindent [D2] P. Deligne, {\it La Conjecture de Weil II}, Publ.
Math. IHES 52 (1980), 137-252.

\bigskip
\noindent [FW] L. Fu and D. Wan, {\it L-functions for symmetric
products of Kloosterman sums}, J. Reine Angew. Math., 589(2005),
79-103.
\bigskip

\noindent [GK] E. Grosse-Kl\"onne, {\it On families of pure slope
L-functions}, Documenta Math., 8(2003), 1-42.

\bigskip

\noindent [K] N. Katz, {\it Gauss Sums, Kloosterman Sums, and
Monodromy Groups}, Princeton University Press, 1988.

\bigskip
\noindent [R] P. Robba, {\it Symmetric powers of $p$-adic Bessel
equation}, J. Reine Angew. Math., 366(1986), 194-220.

\bigskip
\noindent [S] S. Sperber, {\it $p$-Adic hypergeometric functions and
their cohomology}, Duke Math. J. 44 (1977), no. 3, 535-589.

\bigskip
\noindent [W1] D. Wan, {\it Dwork's conjecture on unit root zeta
functions}, Ann. Math., 150(1999), 867-927.

\bigskip
\noindent [W2] D. Wan, {\it Rank one case of Dwork's conjecture}, J.
Amer. Math., Soc., 13(2000), 853-908.

\bigskip
\noindent [W3] D. Wan, {\it Higher rank case of Dwork's conjecture,}
J. Amer. Math., Soc., 13(2000), 807-852.

\end{document}